%17, 21, 25 no such convenion
\input amstex
%ecf55
\input amsppt.sty
\magnification=\magstep1
\advance\vsize-0.5cm\voffset=-1.2cm\advance\hsize1.5cm\hoffset0cm
\NoBlackBoxes
\define\R{{\Bbb R}}\define\Z{{\Bbb Z}}
\def\Aut{\mathop{\fam0 Aut}}
\def\pr{\mathop{\fam0 pr}}
\def\inc{\mathop{\fam0 i}}

\def\forg{\mathop{\fam0 forg}}

\def\con{\mathop{\fam0 con}}
\def\id{\mathop{\fam0 id}}
\def\lk{\mathop{\fam0 lk}}
\def\Int{\mathop{\fam0 Int}}

\def\im{\mathop{\fam0 im}}

\def\t{\widetilde}

\def\Nos{1}

\def\Resp{2}
 
\def\Desp{3}
\def\exli{4\ }
\def\Exli{4}

\def\Desu{5}
\def\gpsi{6\ }
\def\Gpsi{6}

\def\Reso{7}

\def\Depm{8}
\def\sisy{9\ }
\def\Sisy{9}
\def\cosy{10\ }
\def\Cosy{10}
\def\iso{11\ }
\def\Iso{11}

\def\Resy{12}

\def\Deun{13}
\def\adt{14\ }
\def\Adt{14}

\def\Radt{15}
\def\ada{16\ }
\def\Ada{16}
\def\ac{17\ }
\def\Ac{17}

\def\Rmk{18}

\def\Adex{19}
\def\al{20\ }
\def\Al{20}
\def\ut{21\ }
\def\Ut{21}
\def\bast{22\ }
\def\Bast{22}

\def\rl{1.4\ }
\def\Rl{1.4}
\def\is{1.4\ }
\def\Is{1.4}
\def\acl{1.8\ }
\def\Acl{1.8}
\def\su{2.1\ }
\def\Su{2.1}
\def\ind{2.2\ }
\def\Ind{2.2}
\def\bi{2.8\ }
\def\Bi{2.8}

\topmatter
\title How do autodiffeomorphisms act on embeddings?
%An action on the set of
\endtitle
\author A. Skopenkov \endauthor
%$\beta$-invariant, smoothing, self-intersection, almost concordance,
\thanks
This work is supported in part by the Russian Foundation for Basic Research Grants No. 12-01-00748-a and 15-01-06302,
by Simons-IUM Fellowship and by the D. Zimin's Dynasty Foundation Grant.
\newline
I am grateful to S. Avvakumov, D. Crowley, M. Skopenkov and anonymous referee for useful comments.
%\address
\newline
{\it Addresses.} Moscow Institute of Physics and Technology, Faculty of Innovations and High Technology,
Institutskiy per., Dolgoprudnyi, 141700, Russia,
and Independent University of Moscow, B. Vlasy\-ev\-skiy, 11, Moscow, 119002, Russia.
e-mail: skopenko\@mccme.ru
%\endaddress
\endthanks

\abstract
We work in the smooth category.
The following problem was suggested by E. Rees in 2002: describe the precomposition action of self-diffeomorphisms of
$S^p\times S^q$ on the set of isotopy classes of embeddings $S^p\times S^q\to\R^m$.

Let $g:S^p\times S^q\to\R^m$ be an embedding such that
$g|_{a\times S^q}:a\times S^q\to\R^m-g(b\times S^q)$ is null-homotopic for some pair of different points $a,b\in S^p$.

{\bf Theorem.} {\it If $\psi$ is an autodiffeomorphism of $S^p\times S^q$ identical on a neighborhood of
$a\times S^q$ for some $a\in S^p$ and $p\le q$ and $2m\ge 3p+3q+4$, then $g\circ\psi$ is isotopic to $g$.}

Let $N$ be an oriented $(p+q)$-manifold and $f,g$ isotopy classes of embeddings $N\to\R^m$, $S^p\times S^q\to\R^m$, respectively.
As a corollary we obtain that under certain conditions {\it for orientation-preserving embeddings
$s:S^p\times D^q\to N$ the $S^p$-parametric embedded connected sum $f\#_sg$ depends only on
$f,g$ and the homology class of $s|_{S^p\times0}$.}
\endabstract
\endtopmatter

\document
\head 1. Introduction and main results \endhead

{\bf 1.1. Statements of the main results}

This paper is on the classical Knotting Problem: {\it for an $n$-manifold $N$
and a number $m$ describe the set $E^m(N)$ of isotopy classes of embeddings $N\to\R^m$}.
For recent surveys see [Sk08, MA]; whenever possible we refer to these surveys not to original papers.
If the category (PL or smooth) is not mentioned, then the smooth category is tacitly meant.
We denote by $[f]$ the isotopy class of the embedding $f$, except in \S2.1.

An interesting problem is to describe the `precomposition' action of the group $\Aut(N)$ of autodiffeomorphisms of $N$
on $E^m(N)$:
\footnote{
The set of submanifolds of $\R^m$, diffeomorphic to $N$, up to isotopy, is the quotient of $E^m(N)$ by this action.
Action of the group $\Aut_+(N)$ of orientation-preserving autodiffeomorphisms of oriented $N$ is
analogously related to the set of oriented submanifolds of $\R^m$, orientably diffeomorphic to $N$.}
$$\Aut(N)\times E^m(N)\to E^m(N)\quad\text{defined by}\quad (\varphi,[f])\mapsto [f]\circ\varphi:=[f\circ\varphi].$$
E.g. the action of $\Aut(S^1)$ on $E^3(S^1)$ is the same as widely studied `change of the orientation' action.
For the action of $\Aut(S^2)$ on $E^4(S^2)$ see
%papers of Montesinos, Iwase and Hirose, references [5,6,10] of
%Montesinos' is trivial
[Mo83, Hi93, Iw90]. 
For
$$N=T^{p,q}:=S^p\times S^q$$
the problem was raised by E. Rees in 2002.
We obtain a partial solution of this problem.
The main results are Theorems \Exli, \Cosy, \Iso, \Adt, \ada and Corollary \Ac.
The remarks of this text are not used in the statements or proofs of the main results (except that Remarks \Reso.ac
are used for Corollary \ac and Theorem \Exli, respectively).

\smallskip
{\bf Definition \Nos.}
Denote by

$\bullet$ $\sigma:S^q\to S^q$ the reflection w.r.t. the hyperplane $x_1=0$.

$\bullet$ $\inc_q$ the isotopy class of the inclusion $S^q\to S^m$.

$\bullet$ $f\#g$ the embedded connected sum of embeddings or isotopy classes, or the connected sum of autodiffeomorphisms $f$ and $g$.
(See [Av14, \S1] for an accurate definition of the embedded connected sum of embeddings analogous to [Ha66, \S3].
In codimension at least 3 by general position the embedded connected sum of embeddings defines the embedded connected sum of their isotopy classes.)

$\bullet$ $0_k$ the vector of $k$ zero coordinates and $*:=(1,0_{k-1})\in D^k$.

$\bullet$  $\overline\inc=\overline\inc_{p,q,m}:T^{p,q}\to S^m$ {\it the standard embedding} defined by
$\overline\inc(x,y):=(y,0_{m-p-q-1},x)/\sqrt2$, or its abbreviations.

%$$\overline\inc=\overline\inc\phantom{}_{p,q,m}:T^{p,q}\to S^m\quad\text{by}\quad %\overline\inc(x,y):=(y,0_{m-p-q-1},x)/\sqrt2.$$
%Denote by the same notation `$\overline\inc$' abbreviations of $\overline\inc$ and by
$\bullet$ $\inc$ the isotopy class of $\overline\inc$.

\smallskip
{\bf Remark \Resp.}
(a) For the standard embedding $\overline\inc:T^{1,1}\to S^3$ and an autodiffeomorphism $\psi$ of $T^{1,1}$ corresponding to a non-trivial element of $SL_2(\Z)$ we have $\inc\circ\psi\ne \inc$.
This is proved comparing $g|_{*\times S^1}$, $g|_{S^1\times*}$ and
%\linebreak
$\lk(g(*\times S^1),g(-*\times S^1))$ for $g=\inc,\inc\circ\psi$.

For the action of $\Aut_+(T^{1,1})$ on $E^4(T^{1,1})$ see [Hi02]. 
For the case of embeddings $T^{1,2}\to\R^6$, $T^{2,2}\to\R^7$ and $T^{1,3}\to\R^7$ the results of [Sk08', Sk10, CS11, CS] could be useful.
% to attack the Rees problem.
%Cf. [Hi11].

(b) A group structure on $E^m(S^q)$ is defined in [Ha66] for $m\ge q+3$.
%The zero element is represented by $\inc_q$.
The group $E^m(S^q)$ is trivial for $2m\ge3q+4$ [Ha66].

(c) It would be interesing to know how composition with $\sigma$ acts on $E^m(S^q)$.
Composition with $\sigma$ induces an automorphism of $E^m(S^q)$ for $m\ge q+3$
(this is proved analogously to Theorem \iso below).
%(because the reflection w.r.t. the hyperplane $x_2=0$ is isotopic to $\sigma$ and commutes with $\sigma$).
For $m=q+3=7$ the action is identical [Sk10, Symmetry Remark and footnote in \S3].
%(For $m=q+3=6$ analogous idea does not work because althogh the attaching invariant does not change,
%it assumes values in $\pi_3(G_3,SO_3)\cong\Z_2$.)

(d) If $g\in E^m(S^q)$, $\psi\in \Aut_+(S^q)$ and $m\ge q+3$, then $g\circ\psi=g+(\inc_q\circ\psi)$.
The proof is obtained from the proof of Lemma \gpsi below by changing $S^p$ to a point.

(e) Analogously, take any $g\in E^m(T^{p,q})$ and $\psi\in \Aut_+(S^{p+q})$.

If $m\ge p+q+3$, then $g\circ(\id_{T^{p,q}}\#\psi)=g\#(\inc_q\circ\psi)$, where $\id_{T^{p,q}}\#\psi$ is the `connected sum' autodiffeomorphism  of $T^{p,q}$.
(In the notation of [Av14, \S1] the connected sum of the autodiffeomorphisms $a$ and $b$ is the autodiffeomorphism $a\#b: M\#N\to M\#N$ defined by $a\#b|_{M_0} = a|_{M_0}$, $a\#b|_{N_0} = b|_{N_0}$ and $a\#b|_{S^{p-1}\times I} = l_{S^{p-1}\times I}$.)

For $m\ge 2p+q+3$, $g\#(\inc_q\circ\psi)=g$ if and only if $\inc_q\circ\psi=\inc_q$  by
%[Sk15, Theorem 1.1,
%Sk11, Claim 2.4,
[CRS12, Proposition 5.6].

(f) For $q\ge5$ the set of isotopy classes of $\Aut_+(S^q)$ can be identified with the group $\theta_{q+1}$
of homotopy spheres [Sm61, Ce70].
The `composition with $\inc_q$' map $\partial:\theta_{q+1}\to E^m(S^q)$ is a homomorphism appearing in the exact sequence [Ha66, 1.9].
Hence by [Ha66, 1.9] and [Le65, 7.4] we have the following result.

\smallskip
{\bf Proposition \Desp.}
(a) For $m-3=q\in\{7,8,9\}$ there is $\psi\in\Aut_+(S^q)$ such that $\inc_q\circ\psi\ne\inc_q$.
%moreover, $\im\partial$ consists of $2, 2, 3$ elements according to $q=7,8,9$;

(b) For $q\in\{7,8\}$, $m\ge q+4$ and any $\psi\in\Aut_+(S^q)$ we have $\inc_q\circ\psi=\inc_q$.

%Any autodiffeomorphism of $S^n$ is isotopic to an autodiffeomorphism identical on $D^n_+$.
%So the hypothesis of Lemma \gpsi holds automatically for $X$ a point.
%Suppose that $X\in\{D^p,S^p\}$ and  $\cases m\ge p+q+3 & X=D^p \\ m\ge 2p+q+3 & X=S^p \endcases.$

\smallskip
{\bf Theorem \Exli.}
{\it For each integer $l\ge2$ there is $g\in E^{4l+3}(T^{1,2l+1})$ and an autodiffeomorphism $\psi$ of $T^{1,2l+1}$ identical on a neighborhood of $*\times S^{2l+1}$ and such that $g\circ\psi\ne g$ while $\inc\circ\psi=\inc$.

Moreover, we can take as $\psi$ the autodiffeomorphism constructed from the non-trivial element of $\pi_1(SO_{2l+2})$ as in Remark \Reso.c.}

\smallskip
{\bf Definition \Desu.}
An `$S^p$-parametric connected sum' abelian group structure on $E^m(T^{p,q})$ for $m\ge 2p+q+3$ is defined in
[Sk15, \S2.1, MAP] (the definition is recalled in \S2.3).
 
\smallskip
{\bf Lemma \Gpsi.}
{\it If $g\in E^m(T^{p,q})$, $m\ge2p+q+3$ and $\psi$ is an autodiffeomorphism of $T^{p,q}$ identical
on a neighborhood of $S^p\times*$, then $g\circ\psi=g+(\inc\circ\psi)$.}

\smallskip
{\bf Remark \Reso.}
(a) The group $E^m(T^{p,q})$ is trivial for $m\ge p+q+\max\{p,q\}+2$ by the Haefliger Unknotting Theorem
[Sk08, Theorem 2.8.b].
Hence {\it the group structure, Lemmas \Gpsi, \sisy and Theorems \Cosy, \iso below are only interesting for $p\le q$.} So a reader might assume that $p\le q$, although this condititon is not required and so is not added to the
statements of the above-listed results.

The lowest dimensional case when $2p+q+3\le m<p+q+\max\{p,q\}+2$ is $m=8$, $p=1$ and $q=3$.

(b) The condition that $\psi$ is identical on a neighborhood of $S^p\times*$ is essential in Lemma \gpsi by Theorem  \Exli.
%If $\pi_p(SO_{q+1})\ne0$, then the construction of (c) with $p$ and $q$ exchanged apparently gives more   %autodiffeomorphisms of $T^{p,q}$ which are not isotopic to autodiffeomorphisms identical on a neighborhood of %$S^p\times*$.

(c) Let $\varphi:S^q\to SO_{p+1}$ be a smooth map which maps a neighborhood of $*\in S^q$ to the identity.
Define an automorphism
$$\overline\varphi\quad\text{of}\quad T^{p,q}\quad\text{by}\quad \overline\varphi(a,b):=(\varphi(b)a,b).$$
It is clear that $\overline\varphi$ is identical on a neighborhood of $S^p\times*$.

(d) It follows from [MR71, Corollary C] that for $n$ odd and $n\ne3,5,9$, there are pairs of embeddings $f_0,f_1:S^2\times S^{n-2}\to \R^{2n-2}$ such that the normal bundle of $f_0$ is trivial, whereas the normal bundle of $f_1$ is not.
Consequently, in this case the action of $\Aut_+(S^2\times S^{n-2})$ on $E^{2n-2}(S^2\times S^{n-2})$ is not transitive.

\smallskip
{\bf Definition \Depm.}
Define $D^p_+,D^p_-\subset S^p$ by equations $x_1\ge0$, $x_1\le0$, respectively.
Then $S^p=D^p_+\cup D^p_-$.
For an autodiffeomorphism $\alpha$ of $S^q$ denote by
$$\widehat \alpha:=\id\phantom{}_{S^p}\times \alpha$$
the product autodiffeomorphism of $T^{p,q}$.
An autodiffeomorphism $\psi$ of $T^{p,q}$ is {\it symmetric} if
$$\psi(S^p\times D^q_\pm)=S^p\times D^q_\pm\quad\text{and}\quad\widehat\sigma\circ\psi=\psi\circ\widehat\sigma.$$

{\bf Lemma \Sisy.}
{\it If $m\ge2p+q+2$ and $\psi$ is a symmetric autodiffeomorphism of $T^{p,q}$,
then $\inc\circ\psi=\inc$.}

\smallskip
Lemmas \gpsi and \sisy imply the following result.

\smallskip
{\bf Theorem \Cosy.} {\it If $m\ge2p+q+3$ and $\psi$ is a symmetric autodiffeomorphism of $T^{p,q}$ identical
on a neighborhood of $S^p\times*$, then $g\circ\psi=g$ for each $g\in E^m(T^{p,q})$.}

\smallskip
{\bf Theorem \Iso.}
{\it For $m\ge2p+q+3$ the composition with a symmetric autodiffeomorphism of $T^{p,q}$ defines an
automorphism of the group $E^m(T^{p,q})$.}

\smallskip
{\bf Remark \Resy.}
(a) The property of being symmetric depends on the order of factors in $S^p\times S^q$,
i.e. a symmetric autodiffeomorphism of $T^{p,q}$ need not be such of $T^{q,p}$.

(b) For a smooth map $\varphi:S^p\to SO_q\subset SO_{q+1}$ the autodiffeomorphism of $T^{p,q}$ defined analogously to Remark \Reso.c with $p$ and $q$ exchanged, is symmetric and is identical on a neighborhood of $*\times S^q$, but is not necessarily identical on a neighborhood of $S^p\times*$.
%unless  $\varphi$ is homotopy trivial.
%and presumably not isotopic to such
%It would be interesting to know if for $p<q$ any symmetric orientation-preserving autodiffeomorphism of $T^{p,q}$
%is isotopic to a symmetric autodiffeomorphism identical on a neighborhood of $S^p\times*$.

(c) It would be interesting to know if there is a symmetric autodiffeomorphism of $T^{p,q}$ identical
on a neighborhood of $S^p\times*$ but not isotopic to the identity.
If there are none, Theorem \cosy is not interesting.

(d) It would be interesting to know if either of the conditions (that $\psi$ is symmetric or
identical on a neighborhood of $S^p\times*$) is essential in Lemma \sisy and Theorems \Cosy, \Iso.
%i.e. that there are $m,p,q$ and an autodiffeomorphism $\psi$ of $T^{p,q}$ identical
%on a neighborhood of $S^p\times*$ such that $m\ge2p+q+3$ and $\inc\circ\psi\ne\inc$.

(e) Theorem \iso is not covered by Theorem \cosy for any symmetric orientation-reversing autodiffeomorphism of $T^{p,q}$.
%, e.g. for $\widehat\sigma\circ\overline\varphi$, where $\varphi:S^p\to SO_q$ is a map.

(f) The proofs of Lemma \gpsi and Theorems \Cosy, \iso are not hard (\S2).
However, working with distinct embeddings having the same image requires care.
Lemma \sisy is an easy corollary of an unknotting theorem and smoothing theory (\S2).

\smallskip
{\bf Definition \Deun.}
An embedding $g:T^{p,q}\to\R^m$ (or its isotopy class) is called {\it unlinked} if
$$g|_{*\times S^q}:*\times S^q\to\R^m-g((-*)\times S^q)$$
is null-homotopic.
E.g. `$\inc$' is unlinked and
%by general position
any $g\in E^m(T^{p,q})$ is unlinked for $m\ge 2q+2$.
%and for $p=1$, $m=2q+1=4l+1\ge9$ (because ). (Note that any embedding $g:T^{1,2l}\to\R^{4l+1}$ is unlinked.)

\smallskip
{\bf Theorem \Adt.}
{\it Suppose that $p\le q$, \ $g\in E^m(T^{p,q})$ is unlinked and $\psi$ is an autodiffeomorphism of $T^{p,q}$ identical on a neighborhood of $*\times S^q$.

(a) If $2m\ge 3p+3q+4$, then $g\circ\psi=g$.

(b) If $m\ge p+q+2+\max\{p,q/2\}$, then $g\circ\psi=g\#u$ for some $u\in E^m(S^{p+q})$.}

\smallskip
{\bf Remark \Radt.}
(a) Theorem \Adt.a follows from Theorem \Adt.b because $E^m(S^n)=0$ for $2m\ge 3n+4.$

(b) The lowest dimensional case when Theorem \adt holds but is not covered by known results is $m=6$, $p=1$ and $q=2$.
In this case $E^6(S^1\times S^2)$ is described in [Sk08'] and $\Aut_+(S^1\times S^2)$ contains a subgroup $\pi_1(SO_3)$ (cf. Remark \Reso.c).

(c) Theorem \exli falls into the dimension assumption of Theorem \adt and so shows that the unlinkedness assumption in Theorem \adt is essential.

(d) Theorem \Adt.b is obtained from Lemma \al below, an application of the Penrose-Whitehead-Zeeman-Irwin trick
[Sk02, Theorem 2.4] and the Hudson `concordance implies isotopy' result.
See details after Lemma \al below.

(e) The PL analogues of Lemmas \Gpsi, \sisy and Theorems \Cosy, \Iso, \Adt, \ada are true.
The proofs are the same (except that in Lemma \ut we do not need smoothing).

(f) It would be interesting to describe the action of $\Aut(D^p\times S^q)$ on $E^m(D^p\times S^q)$.
% A group structure on $E^m(D^p\times S^q)$ is defined in [Sk15] for $m\ge q+3$.
Analogues for $D^p\times S^q$ of our results on $T^{p,q}$ are correct and could be useful.
%Any or-p autodiffeomorphism of $D^p\times S^q$ is isotopic to an autodiffeomorphism identical
%on a neighborhood of $D^p\times*$?

(g) We do not use any results on $\Aut(T^{p,q})$.
However, the interested reader can find information e.g. in [Le69].
%prop. p. 539 shows that action of Aut_+ = action of \pi_p(SO_{q+1}) for $2m\ge 3n+4$

\bigskip
{\bf 1.2. Application to $S^p$-parametric embedded connected sum}

In the rest of this paper $N$ is a compact $n$-manifold and $f\in E^m(N)$.

%In this paragraph we describe the $S^n$-analogue of our main problem.

For $m\ge n+3$, an embedding $s:D^n\to \Int N$ and $g\in E^m(S^n)$ one can define an embedded connected sum $f\#_sg$
(analogously to the embedded connected sum on $E^m(S^n)$ [Ha66]).
A classical interesting question is

{\it when does $f\#_sg$ depend only on $f,g$ and the component of $N$ containing $s(D^n)$?}
\footnote{
%Clearly, the connected sum depends only on $f,g$ and the isotopy class of $s$.
%Choise of an isotopy class of $s$ is the same as choise of a connected component $C$ of $N$,
%and on an orientation on $C$, if $C$ is orientable.
A sufficient condition for this
%of the orientation
is non-orientability or non-closedness of the component, or $\inc\circ\sigma=\inc$.
We conjecture that this sufficient condition is not necessary.
\newline
Analogous remark should be made for $S^0\times S^n$-analogue of embedded connected sum discussed below.}

If $N$ is connected oriented, then $f\#_sg$ is independent on orientation-preserving $s$ (because every two orientation-preserving embeddings of a disk into a connected oriented manifold are isotopic).
So for $m\ge n+3$ a group structure on $E^m(S^n)$ [Ha66] and an action $\#$ of $E^m(S^n)$ on $E^m(N)$ are defined.
For descriptions of this action see [Sk08', Sk10, CS11, Sk11, CRS12, Proposition 5.6, CS].

%For $m\ge n+3$ and embeddings $g:S^0\times S^n\to\R^m$, $s:S^0\times D^n\to N$ one can define a componentwise
%embedded connected sum $f\#_sg$ (analogously to componentwise embedded connected sum of embeddings
%$S^0\times S^n\to\R^m$ [Ha66]).
%action $\#_s$ of $E^m(S^0\times S^n)$ on $E^m(N)$.
%A classical interesting question is
%{\it when the isotopy class of $f\#_sg$ depends only on $f,g$ and the components of $N$ containing
%$s(S^0\times D^n)$?}

%\footnote{Clearly, the connected sum depends only on the isotopy class of $s$.
%Or, equivalently, on two connected components of $N$, and on an orientation on those of the two which are
%orientable. A sufficient condition for independence of the orientation is triviality of
%the `change of the orientations of components' action on $E^m(S^0\times S^n)$.
%We conjecture that this sufficient  condition is not necessary.}

For $m\ge n+p+3$, an embedding $s:S^p\times D^{n-p}_-\to \Int N$ and $g\in E^m(T^{p,n-p})$ one can define an `$S^p$-parametric embedded connected sum' $f\#_sg$ [MAP, Sk07, pp. 262-264, Sk10', \S2] (see the definition in \S2.3). This defines a group structure on $E^m(T^{p,n-p})$ [Sk06, Sk15] and an action $\#_s$ of $E^m(T^{p,n-p})$ on $E^m(N)$
(because $(f\#_sg)\#_sg'=f\#_s(g+g')$ analogously to [Sk15, \S3, proof of the associativity in the proof of the Group Structure Theorem 2.2]).

We study the following interesting questions (which are classical for $p=0$):

{\it when does $f\#_sg$ depend only on $f,g$ and the isotopy (homotopy, homology) class of $s|_{S^p\times0}$?}

%\footnote{For $p=0$ an isotopy class of $s|_{S^p\times0}$ is the same as a pair of the components of $N$
%together with orientations of those }
%Or the same question for $N$ spin and $s$ preserving the spin structure.\newline
%or only on the homology class $[s]\in H_0(N;\Z_2)$.
% $m=7$ and $N=S^2\times \RP^2\#\C P^2$
%(because `change of the orientation of $S^4$' acts trivially on $E^7(S^4)$
%$\bullet$ $b:H_0(N)\otimes E^m(S^n\sqcup S^n)\to E^m(N)$ defined by
%$b(x\otimes g)f:=f\#_xg$, where the connected sum is applied to $f$-images of
%$N_0$ and the connected component of $N$ corresponding to $N_0,x$.}

A relation of these questions to action by autodiffeomorphisms is as follows.

%Identify $S^p\times D^{n-p}$ and $S^p\times D^{n-p}_-$ by the standard diffeomorphism.

\smallskip
{\bf Theorem \Ada.} {\it Let $N$ be an $n$-manifold, $f\in E^m(N)$, $g\in E^m(T^{p,n-p})$,
\linebreak
$s:S^p\times D^{n-p}_-\to\Int N$ an embedding and $\psi$ a symmetric autodiffeomorphism of $T^{p,n-p}$.
If $m\ge n+p+3$, then $f\#_sg=f\#_{s\circ\psi|_{S^p\times D^{n-p}_-}}(g\circ\psi)$.}

 %The following corollary comprises six statements (Ia), (Ib), (Ic), (IIa), (IIb), (IIc).

\smallskip
{\bf Corollary \Ac.} {\it (I) Let $N$ be an oriented $n$-manifold, $g\in E^m(T^{p,n-p})$ unlinked,
$s:S^p\times D^{n-p}\to \Int N$ an orientation-preserving embedding, $f\in E^m(N)$ and $m\ge n+p+3$, $2m\ge3n+4$.

(a) The sum $f\#_sg$ depends only on $f,g$ and the isotopy class of $s|_{S^p\times0}$.

(b) If $N$ is $(2p+2-n)$-connected, then
%the sum
$f\#_sg$ depends only on $f,g$ and the homotopy class of $s|_{S^p\times0}$.

(c) If $p\ge2$ and $N$ is $(p-1)$-connected, then $f\#_sg$ depends only on $f,g$ and the homology class of $s|_{S^p\times0}$ in $H_p(N;\Z)$.

(II) Modifications of the statements (a), (b) and (c) above hold with `$f\#_sg$' replaced by `the class of $f\#_sg$ in $E^m(N)/\#$' and `$2m\ge3n+4$' replaced by `$2m\ge3n-p+4$'.}

\smallskip
{\bf Remark \Rmk.}
(a) Let $N$ be an oriented $(p-1)$-connected $n$-manifold and $n\ge2p+2$.
Then the Whitney invariant $W:E^{2n-p+1}(N)\to H_p(N;\Z_{(n-p-1)})$ is bijective.
($W(f)$ is a `difference' between $f$ and certain chosen embedding,
see the definition and the statement in [Sk08, \S2, MAW].)
Here $\Z_{(k)}$ is $\Z$ for $k$ even and $\Z_2$ for $k$ odd.
`The parametric connected sum' construction together with explicit construction
$\tau:\Z_{(n-p-1)}\to E^{2n-p+1}(T^{p,n-p})$ of embeddings [Sk08, \S3.4, MAK] give the inverse of $W$.
That is,

{\it for each $f\in E^{2n-p+1}(N)$ an action $\#$ of $H_p(N;\Z_{(n-p-1)})$ on $E^{2n-p+1}(N)$
is well-defined by $\#([s]\otimes r)f:=f\#_s\tau(r)$ and is free and transitive. }

Indeed, by [Sk10', end of \S2]
$$(**)\qquad W(f\#_sg)=W(f)+[s](W(g)\cap[S^p\times*]),\quad\text{where}\quad W(g)\in H_p(T^{p,n-p};\Z_{(n-p-1)}).$$
Hence $W(f\#_ug)=W((f\#_sg)\#_tg)$ when $[u]=[s]+[t]\in H_p(N;\Z_{(n-p-1)})$.
Since $W$ is injective, $f\#_ug=(f\#_sg)\#_tg$.
This and $(f\#_sg)\#_sg'=f\#_s(g+g')$ imply that $\#$ is an action.
By (**), $\#$ is free.
By the injectivity of $W$ and (**), $\#$ is transitive.

%Since $n\ge2p+2$, we have $2n-p+1\ge \frac{3n}2+2\ge n+p+3$.
%Hence by Corollary \ac $f\#_sg$ depends only on $[s]$.

(b) Denote by $E^m_0(T^{p,n-p})$ the subgroup of unlinked embeddings in $E^m(T^{p,n-p})$.
Under the assumptions of Corollary \Ac.c for $2m\ge 3n+4$ a map
$$H_p(N;\Z)\times E^m_0(T^{p,n-p})\times E^m(N)\to E^m(N)\quad\text{is well-defined by}
\quad ([s],g,f)\mapsto f\#_sg.$$
%This map gives an action of $[s]\times E^m_0(T^{p,n-p})$ .
We conjecture that this map gives an action of $H_p(N;E^m_0(T^{p,n-p}))$.

%The formula $([s],f)\mapsto f\#_sg$ does not give an action of $H_1(N;\Z)$ on $E^{2n-1}(N)$ for
%$N=T^{1,n-1}\# T^{1,n-1}$ and certain linked $g\in E^{2n-1}(T^{1,n-1})$, cf. [CS, realization of values
%of $\lambda$]. (why well-defined?) n even

(c) Corollary \ac is trivial for $p=0$.
We conjecture that the assumption $p\ge2$ is superfuous in Corollary \Ac.c.

(d) Fix a certain smooth triangulation of $N$.
Embed $N$ into $\R^M$ for some large $M$.
Denote by $ON$ a tubular neighborhood and by $\nu_N:ON\to N$ the normal bundle of $N$ in $\R^M$.
A {\it stable normal framing} on a subset $X\subset N\subset\R^M$  is an embedding $\zeta:X\times D^{M-n}\to ON$
such that $\zeta(a\times D^{M-n})=\nu_N^{-1}(a)$ for each $a\in X$.
\footnote{A stable normal framing on the $(p+1)$-skeleton of $N$ is the same as
a lifting $N\to BO\left<p+1\right>$ of the stable Gauss map $N\to BO$ [Kr99].
%and the same holds in the PL category replacing $BO$ by $BPL$.
A stable normal framing over the 0-skeleton extendable to the 1-skeleton  is the same as an orientation on $N$.
A stable normal framing over the 1-skeleton extendable to the 2-skeleton is the same as a spin structure on $N$.}
An embedding $s:S^p\times D^{n-p}\to N$ is {\it $\zeta$-good} if $\zeta\circ s$ is the standard stable normal framing of $S^p$.

{\it Let $\zeta$ be a stable normal framing of an open neighborhood $U$ in $N$ of the $(p+1)$-skeleton
(of the triangulation) of $N$.
For $\zeta$-good embeddings $s:S^p\times D^{n-p}\to \Int N$, \ $f\#_sg$ depends only on
$f\in E^m(N),g\in E^m(T^{p,n-p})$ and the isotopy class of $s|_{S^p\times0}$.}
(Then analogously to Corollary \Ac.bc one obtains analogous assertions for homotopy and homology class of $s|_{S^p\times0}$.)

{\it Proof.}
Take two good embeddings $s,s':S^p\times D^{n-p}\to N$ isotopic on $S^p\times0$.
By general position we may assume that the image of the isotopy between $s|_{S^p\times0}$ and $s'|_{S^p\times0}$
is contained in $U$.
Hence the isotopy can be extended to an isotopy of $S^p\times D^{n-p}$ between $s$ and a good embedding
$s'':S^p\times D^{n-p}\to N$ coinciding with $s'$ on $S^p\times0$.
Since both $s'$ and $s''$ are good, they are isotopic.
Hence $s$ and $s'$ are isotopic.
\qed

\head 2. Proofs \endhead

{\bf 2.1. Proof of Theorem \Exli}

Denote by $[a]$ the homotopy class of a map $a$.
The following lemma is possibly known.

\smallskip
{\bf Lemma \Adex.} {\it For each $n\ge2$ there is an autodiffeomorphism $\psi$ of $T^{1,n-1}$ identical on a neighborhood of $*\times S^{n-1}$ for which $\pr_{S^{n-1}}\circ\psi$ is not homotopic to $\pr_{S^{n-1}}$.}

%essentially proved in the proof of the following presumably well-known result This lemma implies that there is
%an autodiffeomorphism of $T^{1,n-1}$  acting trivially on homology but not homotopic to the identity map [Le69].
%{\bf DC: can you give a reference?}}
%\newline
%An autodiffeomorphism of $T^{p,q}$ identical on a neighborhood of $S^p\times*$ induces the identity map in
%homology (this is proved using homology and, for $p=q$, intersection form).
%If $p\le q$, then the converse is apparently false.
%If $p>q$ and $\pi_q(SO_p) = 0$, then the converse is clearly true.}

\smallskip
{\it Proof.} Let $\varphi:S^1\to SO_n$ be a homotopy non-trivial map which maps a neighborhood of $*\in S^1$ to
the identity.
Define an automorphism $\psi$ of $T^{1,n-1}$ by $\psi(a,b):=(a,\varphi(a)b)$.
Clearly, $\psi$ is identical on a neighborhood of $*\times S^{n-1}$.
Let $SG_n$ be the space of maps $S^{n-1}\to S^{n-1}$ of degree +1, the base point being the identity.
%$\id_{S^{n-1}}$.
Identify by the exponential law $\pi_1(SG_n)$ and the set of maps $S^1\times S^{n-1}\to S^{n-1}$
mapping $(*,x)$ to $x$ for each $x\in S^{n-1}$, up to homotopy through such maps.
Let $i:\pi_1(SO_n)\to\pi_1(SG_n)$ the inclusion-induced map.
It is known that $i$ is an isomorphism and $[\pr_{S^{n-1}}\circ\psi]=i[\varphi]\ne i[*]=[\pr_{S^{n-1}}]$.
%$J:\pi_1(SO_n)\to\pi_n(S^{n-1})$
%$\psi$ is not homotopic to the identity map or (equivalently because $\psi$ is identical
%on a neighborhood of $1\times S^{n-1}$) that
So $\pr_{S^{n-1}}\circ\psi$ is not homotopic to $\pr_{S^{n-1}}$.
\qed

%I don't know a simple reference for such a statement.  Viewing the map $\varphi$ as the clutching function
%for an (n-1)-sphere bundle over $S^2$, then it would be enough to show that the total space of this bundle
%is not homotopy equivalent to the total space of the trivial bundle $S^{n-1}\times S^2$.
%This is proven in the attached paper of James and Whitehead - see the table on the final page with n=2
%(in their notation, your n-1 = q).

\smallskip
{\it Construction of example from Theorem \Exli.}
Denote $n:=2l+2$.
Take an autodiffeomorphism $\psi$ of $T^{1,n-1}$ given by Lemma \Adex.
Let $v:S^{n-1}\to S^{n-1}$ be a unit length tangent vector field on $S^{n-1}$ whose degree is +1.
I.e. $v$ is a degree +1 map such that $v(x)\perp x$ for each $x\in S^{n-1}$.
Let $\overline g$ be the composition
$$T^{1,n-1}\overset{\widehat v}\to\to T^{n-1,n-1}\overset{\overline\inc}\to\to \R^{2n-1},\quad\text{where}\quad
\widehat v(e^{i\theta},x):=(v(x)\cos\theta+x\sin\theta,x).$$
Since $v(x)\perp x$ for each $x\in S^{n-1}$, the map $\widehat v$ is well-defined.
Let $g$ be the isotopy class of $\overline g$.

\smallskip
By Theorem \Adt.a $\inc\circ\psi=\inc$.
For the proof that $g\circ\psi\ne g$ we need some preliminaries.

\smallskip
{\it Definition of $L(f)\in\pi_n(S^{n-1})$ for an embedding $f:T^{1,n-1}\to\R^{2n-1}$ coinciding with $\overline g$ on $D^1_+\times S^{n-1}$.}
Denote by $L'(f)$ the  homotopy class of the composition
$$(f\circ(\sigma|_{D^1_+}\times\id\phantom{}_{S^{n-1}}))\cup \overline g|_{D^1_-\times S^{n-1}}:
T^{1,n-1}\to S^{2n-1}-\overline g(1\times S^{n-1})\overset h\to\simeq S^{n-1},$$
where $h$ is a homotopy equivalence of degree $+1$ [Sk08, \S3, MAL].
%Since $\lk(g|_{1\times S^{n-1}},g|_{-1\times S^{n-1}})=\deg v=+1$, we have $L(g)=[\pr_{S^{n-1}}]$. Analogously,
For each $f$ as above the restriction of the above composition to  $*\times S^{n-1}$ is homotopic to $\pr_{S^{n-1}}|_{*\times S^{n-1}}$.
%Hence the composition is homotopic to a map $L'(f):T^{1,n-1}\to S^{n-1}$ whose restriction to $1\times S^{n-1}$
%is $\pr_{S^{n-1}}$.
%Such a map is the same as a map $L'(f):S^1\to SG_n$.
%Invariant $L(f)$ is essentially the same as $[L'(f)]\in\pi_1(SG_n)$.}
%This homotopy class is known to be independent on the choise of $L'(f)$ for fixed $f$.}
Consider the Barrat-Puppe exact sequence:
$$\pi_n(S^{n-1})\overset\#\to\to [T^{1,n-1},S^{n-1}]\overset r\to\to [1\times S^{n-1},S^{n-1}],$$
where $r$ is the restriction and $\#$ extends to the `top cell' action of $\pi_n(S^{n-1})$ on $[T^{1,n-1},S^{n-1}]$.
It is well-known that this action is free, see e.g. [RSS05].
Hence there is a unique class $L(f)\in\pi_n(S^{n-1})$ such that $[L'(f)]=[\pr_{S^{n-1}}]\#L(f)$.

\smallskip
{\it Definition of map $\mu:\pi_{p+q}(S^{m-q-1})\to E^m(T^{p,q})$ for $2m\ge3p+3q+4$.}
For each $x\in \pi_{p+q}(S^{m-q-1})$ take a map $x'$ such that
$$S^{p+q}\overset{x'}\to\to S^m-\overline\inc (D^{p+1}\times S^q)\overset h\to\to S^{m-q-1}$$
represents $x$.
Here $h$ is a homotopy equivalence of degree $+1$ [Sk08, \S3, MAL] and $\inc(x,y):=(y\sqrt{2-|x|^2},0_l,x)/\sqrt2$.
Since $2m\ge3p+3q+4$, there is unique up to isotopy embedding $x''$ homotopic to $x'$.
Let $\mu(x)$ be the isotopy class of $\inc\# x''$.
%here we need p-restriction

\smallskip
{\it Proof that $L(f)$ is an isotopy invariant of $f$ for $n\ge6$ even.}
\footnote{It is not clear that $L(f)$ is preserved through an isotopy of $f$ non-identical on $D^1_+\times S^{n-1}$.
\newline
We conjecture that $L(f)$ is an isotopy invariant of $f$ for $n=4$ (then Theorem \exli holds for $l=1$).
For a proof one possibly needs the results of [CS].
Note that $L(f)=\beta(f)$ for the more complicated $\beta$-invariant of [CS].}
%For each $x\in\pi_n(S^{n-1})$ by definition of sum (recalled in (*)before Lemma \ut below) there is a
%representative $f$ of %the sum of $\mu x$ and the isotopy class of $g$ such that $f=g$ on $D^1_+\times S^{n-1}$
%and  $L(f)=x$. - not suff
For each $x\in\pi_n(S^{n-1})$ by definition of sum (recalled before proof of Lemma \gpsi below) $f$
is a representative of $\mu L(f)+g$.
Hence it suffices to prove that $\mu$ is injective.

% {\bf Lemma \Clt.}
{\it For $n\ge6$ there is the following commutative (up to sign) diagram:
$$\minCDarrowwidth{0pt}\CD
\pi_n(S^{n-1}) @>> \Delta > \pi_{n-1}(S^{n-2}) \\
%@>> \mu'' > \pi_q(V_{m-q,p+1})   \\
@VV \tau V @VV \Sigma V \\
%@VV \tau V  \\
E^{2n}(D^1\times S^n) @>> \lambda > \pi_n(S^{n-1}) @>> \mu > E^{2n-1}(T^{1,n-1})  \endCD.$$
Here $\Delta$ is the map from the exact sequence of the `forgetting the last vector' bundle
$S^{n-2}\to V_{n,2}\to S^{n-1}$, the lower line is exact and $\tau$ is an isomorphism.}
% and $\mu$ is defined below.}

This follows by [Sk06, Lemma 5.1 and Restriction Lemma 5.2] for $p=1$, $q=n-1$ and $m=2n-1$, because
$2(2n-1)\ge 3n+4>10$ for $n\ge6$, so by the smooth version of [Sk02, Theorem 2.4] the map $\forg$ in [Sk06, p.15]
is an isomorphism `respecting' the map $\mu$
%(Cf. [Sk11, Theorem 1.6];
(definitions of $\tau,\lambda$ [Sk06, \S5] are not used here).

We have $\Delta\Sigma x=(1-(-1)^n)x=0$ for each $x\in\pi_{n-1}(S^{n-2})$ [JW54].
Since $n<2(n-1)-2$, the map $\Sigma$ is an isomorphism.
%Identify $\pi_{q+1}(V_{m-q,1})$ with $\pi_{q+1}(S^{m-q-1})$.
Hence $\Delta=0$.
Since both $\tau$ and $\Sigma$ are isomorphisms, this implies that $\lambda=0$.
Hence by exactness $\mu$ is injective.
\qed

\smallskip
{\it Proof that $g\circ\psi\ne g$ for $n\ge6$ even.}
We may assume that $\psi$ is identical on $D^1_+\times S^{n-1}$.
Hence $\overline g\circ\psi=\overline g$ on $D^1_+\times S^{n-1}$.
Thus $L(\overline g\circ\psi)$ is defined.
Clearly, $L'(\overline g\circ\psi)=L'(\overline g)\circ[\psi]=[\pr_{S^{n-1}}\circ\psi]\ne [\pr_{S^{n-1}}]$ by Lemma \Adex.
Hence $L(\overline g\circ\psi)\ne0$.
Thus $g\circ\psi\ne g$.
\qed

%$$[L'(g\circ\psi)]=[L'(g)\circ\psi]=[\pr\phantom{}_{S^{n-1}}\circ\psi]\overset{(*)}\to=
%[\pr\phantom{}_{S^{n-1}}]\#\eta \quad\Rightarrow\quad L(g\circ\psi)=\eta\ne0.$$
%Here $\eta\in\pi_n(S^{n-1})$ is the generator and (*) holds by the above-cited result on the Barrat-Puppe exact %sequence because $\pr_{S^{n-1}}\circ\psi$ is not homotopic to $\pr_{S^{n-1}}$ and $\psi$ is identical on
%$1\times S^{n-1}$.

\bigskip
%\smallskip
{\bf 2.2. Proof of Theorem \Adt.b}

The {\it self-intersection set} of a map $F:N\to\R^m$ is
$$\Sigma(F):=\{x\in N\ |\ \#F^{-1}F(x)\ge1\}.$$

%\smallskip
{\bf Lemma \Al.}
{\it Let $g:T^{p,q}\to\Int B^m$ be an unlinked embedding and $\psi$ an autodiffeomorphism of $T^{p,q}$
identical on a neighborhood of $*\times S^q$.
Then there are a neighborhood $\Delta$ of $*\times S^q$ and a homotopy $G:T^{p,q}\times I\to B^m\times I$ between $g\circ\psi$ and $g$ such that $\Sigma(G)\subset(T^{p,q}-\Delta)\times I$ and $G|_{\Delta\times I}$ is the identical homotopy.}

%For each unlinked embedding $g:T^{p,q}\to\R^m$ the embedding  $g\circ\varphi$ is almost concordant to $g$.

\smallskip
{\it Proof.}
We may assume that $\psi$ is identical on $D^p_+\times S^q$.
The abbreviations
$$g_1,g_2:D^p_-\times S^q\to B^m-g(*\times S^q)\quad\text{of}\quad g \quad\text{and}\quad g\circ\psi$$
coincide on the boundary.
So they form together a map
$$g_{12}:= g\circ(\sigma\times\id\phantom{}_{S^q}\cup\psi|_{D^p_-\times S^q}):T^{p,q}\to B^m-g(*\times S^q).$$
This map factors through the inclusion $g(D^p_-\times S^q)\to B^m-g(*\times S^q)$.
Since $g$ is unlinked, this inclusion is null-homotopic.
Hence there is a null-homotopy $G_{12}$ of $g_{12}$.
Denode by $\con X:=X\times I/X\times 1$ the cone of $X$.
Take the composition
$$D^p_-\times S^q\times I=D^p\times I\times S^q\overset{\alpha\times\id_{S^q}}\to\to
(\con S^p)\times S^q=\frac{T^{p,q}\times I}{\{S^p\times y\times1\}_{y\in S^q}}\overset\beta\to\to$$
$$\overset\beta\to\to\frac{T^{p,q}\times I}{T^{p,q}\times1}=
\con T^{p,q}\overset{G_{12}}\to\to B^m-g(*\times S^q),\quad\text{where}$$
\ \quad
$\bullet$ $\alpha$ is the contraction of $x\times I$ to $[x\times0]\in\con S^p$ for each
$x\in \partial D^p=S^{p-1}$; \ $\alpha$ maps $\partial(D^p\times I)$ to the base $[S^p\times0]$ of the cone and $0\times1/2$ to the vertex $[S^p\times0]$ of the cone;

%$\alpha$ is 1--1 outside $\partial D^p\times I$

$\bullet$ $\beta$ is the contraction of the quotient of $T^{p,q}\times1$ to the vertex of the cone.

There is a neighborhood $\Delta$ of $*$ in $S^p$ such that $g(\Delta\times S^q)\cap G_{12}(\con T^{p,q})=\emptyset$.

Let $D^p_k=D^p_+,D^p_-$ according to $k=0,1$, respectively.
We have for $k=0,1$
$$D^p_-\times S^q\times k=D^p_-\times k\times S^q\overset{\alpha\times\id_{S^q}} \to\to [D^p_k\times0]\times S^q= [D^p_k\times S^q\times0]\overset\beta\to\to[D^p_k\times S^q\times0],$$
$$G_{12}|_{D^p_+\times S^q\times0}=g_1\circ(\sigma\times\id\phantom{}_{S^q})\quad\text{and}
\quad G_{12}|_{D^p_-\times S^q\times0}=g_2.$$
Hence the above composition
$G_{12}\circ\beta\circ(\alpha\times\id_{S^q})$ is a homotopy between $g_1$ and $g_2$ relative to the boundary.
The `union' of this homotopy with the identical homotopy of $D^p_+\times S^q$ is the required homotopy $G$.
\qed

\smallskip
{\it Proof of Theorem \Adt.b.}
Denote a representative of $g$ by the same letter $g$.
Take $\Delta$ and $G$ given by Lemma \Al.
We may assume that $\Delta=D^p_+\times S^q$ and $\Sigma(G)\subset D^p_-\times S^q\times[\frac13,\frac23]$.
Since $m\ge 2p+q+2$, by general position we may assume that
$\Sigma(G)\subset D^p_-\times D^q_-\times[\frac13,\frac23]$, cf. [Sk07, footnote 6].
This means that $G$ is a proper {\it quasi embedding} (see the definition in [Sk02, \S2]).
%; better terminology is {\it almost embedding} and {\it almost concordance} [Sk07, \S2]).
For $p=0$ Theorem \Adt.b is trivial, so we may assume that $m+1\ge p+q+1+3$.
Also $2(m+1)\ge 3(p+q+1)+2-p+1$.
Therefore we can apply [Sk02, Theorem 2.4] to $G$.
We obtain a PL concordance
$$F\quad\text{between}\quad g\quad\text{and}\quad g\circ\psi \quad\text{such that}\quad F=G\quad\text{on}\quad (T^{p,q}-D^p_-\times D^q_-)\times[\frac13,\frac23].$$
Then $F$ is a smooth embedding on this set.
Denote by $u:S^{p+q}\to B^m\times0$ a smooth embedding representing minus the complete obstruction in $E^m(S^{p+q})$ to smoothing $F$ [BH70, Bo71].
Change concordance $F$ by boundary embedded connected sum with the cone (whose vertex is in $B^m\times(0,1)$) over
the  embedding $u$.
The obstruction to smoothing of the new concordance is zero.
Therefore $g\circ\psi$ is smoothly concordant to $g\#u$.
Hence $g\circ\psi$ is smoothly isotopic to $g\#u$ [Hu69].
\qed

%%%\footnote{For $m\ge n-p+2+\max\{p,n-p\}$ we have
%%%$[S^p\times S^{n-p},S^{m-n+p-1}]\cong\pi_{n-p}(S^{m-n+p-1})\oplus\pi_n(S^{m-n+p-1})$.
%%%For general $g\in E^m(T^{p,n-p})$ we conjecture that $[\Psi]$ goes to
%%%$\lambda(g)\oplus(\lambda(g)\circ j(\varphi))$ under this isomorphism.
%%%Here $\lambda(g)\in\pi_{n-p}(S^{m-n+p-1})$ is the linking coefficient of $g|_{a\times S^{n-p}}$ and
%%%$g(b\times S^{n-p})$ [Sk08, \S3], and $j:\pi_p(SO_{n-p})\to\pi_n(S^{n-p})$ is the map defined by letting $j(x)$
%%%be the framed submanifold formed by the standard $S^p\subset S^n$ and the framing on $S^p$ twisted by $x$.}

\bigskip
%\newpage
{\bf 2.3. Proof of Lemmas \Gpsi, \Sisy, Theorem \iso and Corollary \Ac}

In this subsection we omit the composition sign $f\circ g$, writing $fg$ for $f\circ g$.
%and denote by $\overline g$ a representative of the isotopy class $g$.

Let $N$ be an $n$-manifold and $s:S^p\times D^{n-p}_-\to N$ an embedding.
A map $f:N\to S^m$ is called {\it $s$-standardized} if
$$f(N-\im s)\subset\Int D^m_+\quad\text{and}\quad fs=\overline\inc.$$
Denote by $i:S^p\times D^{n-p}_-\to T^{p,n-p}$
%(the composition of the standard diffeomorphism and)
the inclusion.
Unless for the proof of Theorem \ada at the end of this subsection a reader may assume that $N=T^{p,n-p}$ and $s=i$.

Denote by $R$ the reflection of $\R^m$ with respect to $0\times0\times\R^{m-2}$,
%x_1=x_2=0
and also abbreviations of this reflection.
Since $m\ge n+p+2$, analogously to [Sk15, Standardization Lemma 2.1] there are an $s$-standardized representative $\overline f$ of $f\in E^m(N)$ and an $i$-standardized representative $\overline g$ of $g\in E^m(T^{p,n-p})$.
Then a representative $\overline h$ of $f\#_sg$ is defined by
$$(*)\qquad \overline h(a):=\cases \overline f(a) &a\not\in\im s\\
R\overline g\widehat Rs^{-1}(a)           &a\in\im s\endcases.$$
The two formulas agree on $\partial\im s$ because $\overline\inc=R\overline\inc\widehat R$;
clearly, $h$ is a ($C^1$-smooth) embedding.
For $m\ge n+p+3$

$\bullet$ this gives a well-defined map $E^m(N)\times E^m(T^{p,n-p})\to E^m(N)$
(analogously to [Sk15, \S3, beginning of the proof of the Group Structure Theorem 2.2]).
\footnote{In [Sk10, Definition 1.4 of the action $b$] it was essentially constructed an action
$b:H_p(N;\pi^S_{2n-p-1-m})\to E^m(N)$ for a closed orientable $(p-1)$-connected $n$-manifold $N$ and $2m\ge3n+4-p$.
%(without the $p$-parallelizability condition).
There is a map $\pi^S_{2n-p-1-m}=\pi_{n-p-1}(S^{m-n})\overset\mu\to\to \pi_{n-p-1}(V_{m-n+p,p+1})\overset\tau\to\to E^m(T^{p,n-p})$ [Sk08, 3.4, MAK] whose image consists of unlinked embeddings.
We have $b([s]\otimes x,f)=f\#_s\tau\mu(x)$.
% for each $x\in \pi^S_{2n-p-1-m}$.
The set of unlinked embeddings is $\im(\tau\mu)$ if either $p=1$ or $m\ge2n-p$.
So the action $b$ is the `top cell part' of the map $([s],g,f)\mapsto f\#_s g$.
The `top cell part' is the same as the whole map if either $p=1$ or $m\ge2n-p$.}

$\bullet$ $+:=\#_i$ gives a well-defined abelian group structure on $E^m(T^{p,n-p})$ [Sk15, \S2.1].

Let $R^t$ be the rotation of $\R^m$ whose restriction to the plane $\R^2\times0$ is the rotation through the angle
$+\pi t$ and which leaves the orthogonal complement $0\times\R^{m-2}$ fixed.

%A map $f:X\times S^q\to S^m$ is called {\it standardized} if
%$$f(X\times\Int D^q_+)\subset\Int D^m_+\quad\text{and}\quad f_{X\times\Int D^q_-}=g_-.$$
%For standardized embeddings $f,g:X\times S^q\to S^m$ a representative $h$ of $[f]+[g]$ is defined by
%the formula (*) in which $s=i$.
%A map $g:T^{p,n-p}\to S^m$ is called {\it standardized} if $g(S^p\times\Int D^{n-p}_+)\subset\Int D^m_+$ and
%$g|_{S^p\times D^{n-p}_-}=g_-$.

\smallskip
{\it Proof of Lemma \Gpsi.}
Let $\inc_\psi:=R\overline\inc \psi\widehat R$.
We may assume that $\psi$ is identical on $S^p\times D^q_+$.
This,
$$R(D^m_\pm)=D^m_\mp,\quad\widehat R(S^p\times D^q_\pm)=S^p\times D^q_\mp\quad\text{and}\quad R\overline\inc\widehat R=\overline\inc$$
imply that embedding $\inc_\psi$ is $i$-standardized.
There is an isotopy $R^t\overline\inc\psi\widehat{R^t}$ between $\inc_\psi$ and $\overline\inc\psi$.
By the Standardization Lemma of [Sk07, Sk10'] there is a standardized representative $\overline g$ of $g$.
Then a representative $h$ of $g+\inc\psi$ is defined by (*) for $\overline f,\overline g,s$ replaced by $\overline g,\inc_\psi,i$, respectively.
We have $h=\overline g=\overline g\psi$ on $S^p\times D^q_+$ and
$h=Ri_\psi\widehat R=\overline\inc\psi=\overline g\psi$ on $S^p\times D^q_-$.
Hence $h=\overline g\psi$.
\qed

%An embedding $g:X\times D^{n-p}\to B^m$ is {\it $X$-proper} if $g(X\times\partial D^{n-p})\subset\partial B^m$.

\smallskip
{\bf
%Unknotting
Lemma \Ut.}
{\it Any two proper embeddings $S^p\times D^q\to B^m$ are properly isotopic for $m\ge 2p+q+2$.}
%{\it Let $N$ be a compact $n$-dimensional manifold and $f,g:N\to D^m$ proper embeddings.
%If $m\ge n+3$ and $(N,\partial N)$ is $(2n-m+1)$-connected, then $f$ and $g$ are properly isotopic.}

\smallskip
{\it Proof.} The pair $(S^p\times D^q,S^p\times\partial D^q)$ is $(q-1)$-connected.
Since $m\ge2p+q+2$, this pair is $(2(p+q)-m+1)$-connected.
Therefore any two proper embeddings $S^p\times D^q\to B^m$ are properly PL isotopic [Hu69, Theorem 10.2].
%for $X=D^p$ and $p\ge1$ we need a version!

Obstructions to smoothing this isotopy (moving $S^p\times\partial D^q$ in $\partial B^m$) are in
\linebreak
$H^j(S^p\times D^q;E^{m-p-q+j}(S^j))$ [Ha67, BH70, Bo71].
The only non-trivial obstruction could appear for
%$X=S^p$ and
$j=p$.
Since $m-p-q\ge2p+2$, we have $2(m-q)\ge2(2p+2)\ge3p+4$, so this obstruction is zero.
\qed

\smallskip
In this subsection we denote by the same letter a symmetric autodiffeomorphism of $T^{p,q}$
and its abbreviation $S^p\times D^q_\pm\to S^p\times D^q_\pm$.

Lemma \sisy is implied by  the following Lemma \Sisy'.a.
\footnote{Lemma \sisy is used in the proof of Theorem \Iso.
So although Lemma \sisy for $m\ge 2p+q+3$ follows from Theorem \Iso, Lemma \sisy for $m\ge 2p+q+3$ is not a corollary of Theorem \Iso.}

\smallskip
{\bf Lemma \Sisy'.} {\it Let $\psi$ be a symmetric autodiffeomorphism of $T^{p,q}$.

(a) If $m\ge2p+q+2$, then there is an isotopy $H_t:S^m\to S^m$ of the identity map $H_0$ such that
$H_1\overline\inc\psi=\overline\inc$ and $H_1(D^m_\pm)=D^m_\pm$.
%, (no need?) $H_1\sigma=\sigma H_1$

(b) If $H_t$ is an isotopy from (a), embedding $\overline g:T^{p,q}\to S^m$ is $i$-standardized and embedding
$\overline f:N\to S^m$ is $s$-standardized, then

$\bullet$ embedding $H_1\overline f$ is $s\psi$-standardized;

$\bullet$ embedding $H_1\overline g\psi$ is $i$-standardized;

$\bullet$ embedding $\overline g_\psi:=RH_1R\overline g\widehat R\psi\widehat R$ is $i$-standardized and isotopic to
$\overline g\psi$.}

\smallskip
{\it Proof of (a).} By Lemma \ut there is an isotopy between
$\overline\inc:S^p\times D^q_+\to D^m_+$ and $\overline\inc\psi:S^p\times D^q_+\to D^m_+$.
Since $\psi$ is symmetric and being proper includes being orthogonal near the boundary,
the symmetric extension of the above isotopy w.r.t. the hyperplane $x_1=0$ is a smooth isotopy.
This extension is as required.
\qed

\smallskip
{\it Proof of (b).}
We have
$$H_1\overline fs\psi=H_1\overline\inc\psi=\overline\inc\quad\text{on}\quad S^p\times D^q_-\quad\text{and}\quad
H_1\overline f(N-\im s)\subset H_1(D^m_+)=D^m_+.$$
Thus $H_1\overline f$ is $s\psi$-standardized.

%We have
%$$H_1g\psi(S^p\times D^q_+)=H_1g(S^p\times D^q_+)\subset H_1(D^m_+)=D^m_+\quad\text{and}$$
%$$\text{on}\quad S^p\times D^q_-\quad\text{we have}\quad H_1g\psi=H_1\inc\psi=\inc.$$
%Thus $H_1g\psi$ is $i$-standardized.

Clearly, {\it an embedding $\overline g:T^{p,n-p}\to\R^m$ is $i$-standardized if and only if $\overline g\psi$ is $\psi^{-1}i$-standardized.}
So the second bullet point follows from the first one.

Both $\psi$ and $\widehat R$ preserve $S^p\times D^q_\pm$, both $H_1$ and $R$ preserve $D^m_\pm$, both
$\overline g$ and $H_1\overline g\psi$ are $i$-standardized, and
$\overline\inc=R\overline\inc\widehat R=H_1\overline\inc\psi$.
Hence embedding $\overline g_\psi=R(H_1(R\overline g\widehat R)\psi)\widehat R$ is $i$-standardized.
Then $R^tH_tR^t\overline g\widehat R^t\psi\widehat R^t$ is an isotopy between $\overline g_\psi$ and
$\overline g\psi$.
\qed

%So if embedding $f:N\to S^m$ is $s$-standardized, then
%$$H_1fs\psi(S^p\times D^{n-p}_+)=H_1fs(S^p\times D^{n-p}_+)\subset H_1(D^m_+)=D^m_+\quad\text{and}$$
%$$\text{on}\quad S^p\times D^{n-p}_-\quad\text{we have}\quad H_1fs\psi=H_1g_-\psi=g_-\psi^{-1}\psi=g_-.\quad\qed$$

\smallskip
{\it Proof of Theorem \Iso.}
Precomposition with $\psi$ defines a self-bijection of $E^m(T^{p,q})$.

Take $i$-standardized representatives $\overline g,\overline g'$ of $g,g'\in E^m(T^{p,q})$.
Then a representative $\overline h$ of $g'+g$ is defined by (*) for $\overline f=\overline g'$,
$\overline g=\overline g$ and $s=i$.
Take an an isotopy $H_t$ given by Lemma \Sisy'.a.
Then by Lemma \Sisy'.b $H_1\overline g'\psi$ and $\overline g_\psi$ are $i$-standardized.
Hence by Lemma \Sisy'.b a representative $\overline h_\psi$ of $g'\psi+g\psi$ is defined by (*) for
$\overline f=H_1\overline g'\psi$, $\overline g=\overline g_\psi$ and $s=i$.
We have $g'\psi+g\psi=(g'+g)\psi$ because
$$\cases \overline h_\psi=H_1\overline g'\psi=H_1\overline h\psi &\text{on } S^p\times D^q_+\\
\overline h_\psi=R\overline g_\psi\widehat R=H_1R\overline g\widehat R\psi=H_1\overline h\psi
&\text{on } S^p\times D^q_-\endcases .$$
Thus composition with $\psi$ defines an automorpism of $E^m(T^{p,q})$.
\qed

\smallskip
{\it Proof of Theorem \Ada.}
By the Standardization Lemma of [Sk07, Sk10'] there are $s$-standardized and $i$-standardized representatives  $\overline f$ of $f$ and $\overline g$ of $g$, respectively.
%In this proof we denote by $f,g$ these embeddings.
Then a representative $\overline h$ of $f\#_sg$ is defined by (*).
Take an isotopy $H_t$ given by Lemma \Sisy'.a.
Then by Lemma \Sisy'.b  $H_1\overline f$ is $s\psi$-standardized and $\overline g_\psi$ is $i$-standardized.
Hence by Lemma \Sisy'.b a representative $\overline h_\psi$ of $f\#_{s\psi}g\psi$ is defined by  (*) for $\overline f$, $\overline g$ and $s$ replaced by $H_1\overline f$, $\overline g_\psi$ and $s\psi$, respectively.
We have $f\#_s g=f\#_{s\psi} g\psi$ because
$$\cases \overline h_\psi=H_1\overline f=H_1\overline h & \text{on } N-\Int\im s \\
\overline h_\psi=R\overline g_\psi\widehat R\psi^{-1}s^{-1}=H_1R\overline g\widehat Rs^{-1}=H_1\overline h
& \text{on }\im s
\endcases.\quad\qed$$

%\smallskip
{\it Proof of Corollary \Ac.}
If $2p+2>n$, then $m\ge n+p+3\ge n+\max\{p,n-p\}+2$.
Hence $g=\inc$ by Remark \Reso.a.
Then $f\#_sg=f$.
So we may assume that $2p+2\le n$.

Part (a) follows by Theorems \adt and \ada because $p\le n-p$ and any orientation-preserving embedding
$s:S^p\times D^{n-p}\to N$ extending given $s|_{S^p\times0}$ is isotopic to the composition of one fixed such embedding with an autodiffeomorphism of $S^p\times D^{n-p}$ defined by
$(a,b)\mapsto(a,\varphi(a)b)$ for certain map $\varphi:S^p\to SO_{n-p}$.
\footnote{If either $2\le p\ge n-2$ or ($n\ge2p+2$ and $p\equiv2,4,5,6\mod8$), then $\pi_p(SO_{n-p})=0$.
If $m\ge2n+2-p$ and $n\ge2p$, then $E^m(T^{p,n-p})=0$ by Remark \Reso.a.
The dimension restriction in Corollary \ac can be replaced by any condition of this footnote.}

%Instead of application of Theorem \adt one can imitate its proof!!!

Part (b) follows by (a) and the analogue of the Haefliger Unknotting Theorem [Sk08, Theorem 2.8.b] for embeddings $S^p\to N$ because $2n\ge4p+4\ge3p+4$.

Part (c) follows by (b) and the Hurewicz Theorem because $n\ge 2p+2\ge6 \Rightarrow p-1\ge 2p+2-n$.
\qed

% Note that the proof of the completeness [Sk07] is incorrect (because of incorrect
% reference to [Hu70, Ha84]: there are relative invariant to absolute homotopy there, but not an absolute
% invariant to relative homotopy), but can be recovered using the same ideas [Sk06].

%\newpage


\Refs
\widestnumber\key{CRS0}

\ref \key Av 14 \by S. Avvakumov
\paper The classification of certain linked 3-manifolds in 6-space
\yr \vol  \jour \moreref 	arXiv:1408.3918 [math.GT]
\endref


\ref \key BH70 \by J. Bo\'echat and A. Haefliger \pages 156--166
\paper Plongements diff\'erentiables de vari\'et\'es orient\'ees de dimension 4 dans $\R^7$
\yr 1970 \vol  \jour Essays on topology and related topics (Springer,1970)
\endref

\ref \key Bo71 \by J. Boechat \pages 141--161
\paper Plongements differentiables de varietes de dimension $4k$ dans $\R^{6k+1}$
\yr 1971 \vol 46:2 \jour Comment. Math. Helv.
\endref


\ref \key Ce70 \by J. Cerf
\paper La stratification naturelle des espaces de fonctions diff\'erentiables r\'eelles et le th\'eor\`eme de la pseudo-isotopie
\yr 1970 \vol 39 \jour Inst. Hautes \'Etudes Sci. Publ. Math. \pages 5--173
\endref


%\ref  \key CRS07 \by M. Cencelj, D. Repov\v s and M. Skopenkov
%\paper Homotopy type of the complement of an immersion and classification of embeddings of tori
%\jour Uspekhi Mat. Nauk \vol 62:5 \yr 2007  \pages 165-166
%\moreref English transl: Russian Math. Surveys \vol 62:5 \yr 2007  \pages
%\moreref arxiv:math/0803.4285
%\endref

\ref  \key CRS12 \by M. Cencelj, D. Repov\v s and M. Skopenkov
\paper  Classification of knotted tori in the 2-metastable dimension
\jour Mat. Sbornik  \vol 203:11 \yr 2012  \pages 1654-1681
\moreref arXiv:math/0811.2745
\endref


\ref  \key CS11 \by D. Crowley and A.  Skopenkov
\paper A classification of smooth embeddings of 4-manifolds in 7-space, II
\jour Internat. J. Math. \vol 22:6 \yr 2011 \pages 731-757
\moreref arxiv:math/0808.1795
\endref

\ref  \key CS \by D. Crowley and A.  Skopenkov
\paper A classification of smooth embeddings of 4-manifolds in 7-space, III, draft
\jour  \vol \yr  \pages \moreref
\endref


\ref \key Ha66 \by A. Haefliger
\paper Differentiable embeddings of $S^n$ in $S^{n+q}$ for $q>2$
\pages 402--436 \jour Ann. Math., Ser.3 \vol 83 \yr 1966
\endref

%\ref \key Ha66' \by A.~Haefliger
%\paper Enlacements de spheres en codimension superiure a 2
%\jour Comment. Math. Helv. \vol 41 \yr 1966-67 \pages 51--72
%\endref

\ref \key Ha67 \by A. Haefliger \pages 221--240
\paper Lissage des immersions-I
\yr 1967 \vol 6 \jour Topology
\endref

%\ref \key HH63 \by A. Haefliger and M. W. Hirsch \pages 129--135
%\paper On existence and classification of differential embeddings
%\yr 1963 \vol 2 \jour Topology
%\endref

\ref \key Hi93 \by S. Hirose
\paper On diffeomorphisms over $T^2$-knot
\jour Proc. of A.M.S. \vol 119 \yr 1993 \pages 1009–1018
\endref



\ref \key Hi02 \by S. Hirose
\paper On diffeomorphisms over surfaces trivially embedded in the 4-sphere
\jour Algebr. Geom. Topol.  \vol 2 \yr 2002 \pages 791-824
\moreref arxiv:math/0211019
\endref



%\ref \key Hi11 \by S. Hirose
%\paper On diffeomorphisms over non-orientable surfaces standardly embedded in the 4-sphere
%\jour  \vol  \yr  \pages arxiv:math/1109.1668
%\endref

\ref \key Hu69 \by J. F. P. Hudson \book Piecewise-Linear Topology
\bookinfo \publ Benjamin \publaddr New York, Amsterdam \yr 1969
\endref

\ref \key Iw90 \by Z. Iwase
\paper Dehn surgery along a torus $T^2$-knot, II.
\jour Japan. J. Math. \vol 16 \yr 1990 \pages  171–196
\endref


\ref \key JW54 \by I. M. James and J. H. C. Whitehead
\paper The homotopy theory of sphere bundles over spheres,
\jour Proc. London Math. Soc. (3)  I: \yr 1954 \vol 4 \pages 196--218
\moreref II: \yr 1955 \vol 5 \pages 148--166
\endref

\ref \key Kr99 \by M. Kreck
\paper Surgery and duality
\pages 707--754 \jour Ann. Math. \vol 149 \yr 1999
\endref

\ref \key Le65 \by J.~Levine
\paper A classification of differentiable knots
\jour Ann. of Math.  \vol 82 \yr 1965 \pages 15--50
\endref


\ref \key Le69 \by J. Levine
\paper Self-equivalences of $S^n\times S^k$
\pages 523--543 \jour Trans. Amer. Math. Soc. \vol 143 \yr 1969
\endref

\ref \key MA \by
\paper http://www.map.mpim-bonn.mpg.de/High\_codimension\_embeddings
\pages \jour Manifold Atlas Project (unrefereed page)
\vol \yr
\endref

\ref \key MAK \by
\paper http://www.map.mpim-bonn.mpg.de/Knotted\_tori\#Examples
\pages \jour Manifold Atlas Project (unrefereed page)
\vol \yr
\endref

\ref \key MAL \by
\paper http://www.map.mpim-bonn.mpg.de/Links,\_i.e.\_embeddings\_of\_non-connected\_manifolds\#
\linebreak
General\_position\_and\_the\_Hopf\_linking
\pages \jour Manifold Atlas Project (unrefereed page)
\vol \yr
\endref

\ref \key MAP \by
\paper http://www.map.him.uni-bonn.de/Parametric\_connected\_sum.
\pages \jour Manifold Atlas Project (unrefereed page)
\vol \yr
\endref



\ref \key MAW \by
\paper http://www.map.mpim-bonn.mpg.de/Embeddings\_just\_below\_the\_stable\_range:\_classification\#
\linebreak
The\_Whitney\_invariant
\pages \jour Manifold Atlas Project (unrefereed page)
\vol \yr
\endref

\ref  \key Mo83 \by J. M. Montesinos
\paper On twins in the four-sphere, I
\jour Quart. J. Math Oxford (2) \vol 34  \yr 1983  \pages 171-199
\endref



\ref  \key MR71 \by R. J. Milgram and E. Rees
\paper On the normal bundle to an embedding
\jour Topology \vol 10  \pages 299--308  \yr 1971
\endref


\ref  \key RSS05 \by D. Repov\v s,  M. Skopenkov and F. Spaggiari
\paper   On the Pontryagin-Steenrod-Wu theorem
\jour Israel J. Math.  \vol 145 \yr 2005  \pages 341-347
\moreref arXiv:math/0808.1209
\endref

\ref \key Sk02 \by A. Skopenkov
\paper On the Haefliger-Hirsch-Wu invariants for embeddings and immersions
\yr 2002 \vol \jour Comment. Math. Helv. \pages 78--124
\endref

\ref \key Sk06 \by A. Skopenkov
\paper  Classification of embeddings below the metastable dimension
\yr \vol \jour \pages
\moreref arxiv:math/0607422
\endref

\ref \key Sk07 \by A. Skopenkov
\paper A new invariant and parametric connected sum of embeddings
\yr 2007 \vol 197 \jour Fund. Math. \pages 253--269
\moreref arxiv:math/0509621
\endref

\ref \key Sk08 \by A. Skopenkov
\paper Embedding and knotting of manifolds in Euclidean spaces,
in: Surveys in Contemporary Mathematics, Ed. N. Young and Y. Choi
\yr 2008 \vol 347 \jour London Math. Soc. Lect. Notes \pages 248--342
\moreref arxiv:math/0604045
\endref

\ref \key Sk08' \by  A. Skopenkov
\paper Classification of smooth embeddings of 3-manifolds in the 6-space
\yr 2008 \vol 260:3 \jour Math. Zeitschrift \pages 647-672
\moreref  arxiv:math/0603429, DOI: 10.1007/s00209-007-0294-1
\endref

\ref  \key Sk10 \by A.  Skopenkov
\paper A classification of smooth embeddings of 4-manifolds in 7-space, I
\jour Topol. Appl. \vol 157 \yr 2010 \pages 2094-2110
\moreref arxiv:math/0512594
%http://dx.doi.org/10.1016/j.topol.2010.05.003
\endref


\ref  \key Sk10' \by A.  Skopenkov
\paper  Embeddings of $k$-connected $n$-manifolds into $\R^{2n-k-1}$
\jour Proc. Amer. Math. Soc. \vol 138 \yr 2010 \pages 3377--3389
%\linebreak
\moreref arxiv:0812.0263 [math.GT]
\endref

\ref  \key Sk11 \by M.  Skopenkov
\paper When is the set of embeddings finite?
\moreref arxiv:math/1106.1878
\endref

\ref \key Sk15 \by A. Skopenkov
\paper  Classification of knotted tori
\yr \vol \jour \pages
\moreref arXiv:1502.04470 [math.GT]
\endref


\ref \key Sm61 \by S.~Smale
\paper Generalized Poincare's conjecture in dimensions greater than 4
\yr 1961 \vol 74 \jour Ann. of Math. (2) \pages 391--466
\endref

\endRefs
\enddocument